\documentclass[reqno,12pt,letterpaper]{amsart}
\usepackage{amsmath,amssymb,amsthm,graphicx,mathrsfs,url}
\usepackage[usenames,dvipsnames]{color}
\usepackage[colorlinks=true,linkcolor=Red,citecolor=Green]{hyperref}

\usepackage{tikz}

\usepackage{graphicx}
\usepackage{caption}
\usepackage{subcaption}
\usepackage[colorlinks=true]{hyperref}
\usepackage{color}
\usepackage{amsmath,amsfonts}
\usepackage{calrsfs}
\usepackage{amsmath,environ}
\DeclareMathAlphabet{\pazocal}{OMS}{zplm}{m}{n}

\numberwithin{equation}{section}
\usepackage{amsmath, amsthm}
    \usepackage{dsfont}
    \usepackage{fancybox}
    \usepackage{amssymb}
    \graphicspath{ {images/} }

 	\definecolor{cadmiumgreen}{rgb}{0.0, 0.42, 0.24}

\newcommand{\R}{\mathbb{R}}

\newcommand{\w}{\omega}
\newcommand{\dd}[2]{\dfrac{\partial #1}{\partial #2}}
\newcommand{\matrice}[1]{\left( \begin{matrix}
#1
\end{matrix} \right)}

\newcommand{\1}{\mathds{1}}
\newcommand{\al}{\al}

 \newcommand{\Id}{\mathrm{Id}}
 \newcommand{\tRR}{\widetilde{\RRR}}

\newcommand{\lr}[1]{\langle #1 \rangle}

\newcommand{\Det}{{\operatorname{Det}}}

\newcommand{\GG}{\mathcal{G}}

\newcommand{\mi}{\setminus}
\newcommand{\supp}{\mathrm{supp}}

\newcommand{\Bb}{\mathbb{B}}

\newcommand{\Ss}{\mathbb{S}}
\newcommand{\BB}{\mathcal{B}}

\newcommand{\EE}{\mathcal{E}}

\newcommand{\MM}{\mathcal{M}}
\newcommand{\OO}{\mathcal{O}}

\newcommand{\SSS}{{\mathcal{S}}}
\newcommand{\Hh}{{\mathbb{H}}}

\newcommand{\az}{\alpha}

\newcommand{\Mm}{\mathbb{M}}

\newcommand{\Aff}{\operatorname{Aff}}
\newcommand{\tp}{{\tilde{p}}}

\newcommand{\RRR}{\mathcal{R}}
\newcommand{\RR}{\mathcal{R}}
\newcommand{\tlambda}{{\tilde{\lambda}}}
\newcommand{\de}{\mathrel{\stackrel{\makebox[0pt]{\mbox{\normalfont\tiny def}}}{=}}}
\NewEnviron{equations}{%
\begin{equation}\begin{gathered}
  \BODY
\end{gathered}\end{equation}
}
\newcommand{\arxiv}[1]{\href{http://arxiv.org/abs/#1}{arXiv:#1}}

\NewEnviron{equations*}{%
\begin{equation*}\begin{gathered}
  \BODY
\end{gathered}\end{equation*}
}

\title[Extremizers for $k$-plane transform inequalities.]{Existence and non-existence of extremizers for $k$-plane transform inequalities.}

\author{Alexis Drouot}
\email{alexis.drouot@gmail.com}

\newtheorem{thm}{Theorem}
\newtheorem{proposition}{Proposition}

\newtheorem{lem}{Lemma}[section]

\newtheorem{theorem}[thm]{Theorem}

\begin{document}

\maketitle

\begin{abstract} We provide sharp forms of $k$-plane transform inequalities on the $d$-dimensional sphere $\Ss^d$ and the $d$-dimensional hyperbolic space $\Hh^d$. In particular, we prove that extremizers do not exist for $\Hh^d$. This work is a natural extension of previous results for the $k$-plane transform on $\R^d$.
\end{abstract}

\section{Introduction}

In this note, we study sharp inequalities associated with the $k$-plane transform on the sphere $\Ss^d$ and on the hyperbolic plane $\Hh^d$, defined as follows. Let $\Mm^d_\nu, \ \nu \in \{0,+,-\}$ be the Riemannian manifold equal to $\R^d$ for $\nu=0$, the sphere $\Ss^d$ for $\nu=+$, and the hyperbolic space $\Hh^d$ for $\nu=-$. The spaces $\Mm_\pm^d$ can be seen as subspaces of $\R^{d+1}$:
\begin{equation*}
\Mm_\pm^d = \{ (\zeta',\zeta_{d+1}) \in \R^d \times \R, \ |\zeta_{d+1}|^2 \pm |\zeta'|^2=1 \}. 
\end{equation*}
A $k$-plane in $\Mm_\pm^d$ is defined as the intersection of a $k+1$ plane through the origin in $\R^{d+1}$ with $\Mm_\pm^d$. Let $\MM_k(\Mm_\nu^d)$ be the set of $k$-planes in  $\Mm_\nu^d$. The $k$-plane transform on $\Mm_\nu^d$ is the operator given by
\begin{equation*}
\RR_\nu f(\pi) \de \int_\pi f d\lambda_\pi, \ \ \ \ f \in C_0^\infty(\Mm^d_\nu), \ \ \ \ \pi \in \MM_\nu(\Mm_\nu^d),
\end{equation*}
where $d\lambda_\pi$ is the measure corresponding to the induced Riemannian metric on $\pi$.

The $k$-plane transform on non-flat manifolds was introduced in Helgason \cite{Helgason}. We refer to Berenstein--Casadio--Kurusa \cite{BerTarKur} and Quinto \cite{Quinto} for range characterization and support theorems, Rubin \cite{Rubin} for pointwise inversion formula in the hyperbolic case, and Berenstein--Rubin \cite{BerRub} for Radon transforms of functions defined almost everywhere. In this paper, we focus on sharp $L^p(\Mm_\nu^d)$ to $L^q(\MM_k(\Mm_\nu^d))$ inequalities. Modulo a multiplicative constant, there exists a unique measure on $\MM_k(\Mm_\nu^d)$ that is invariant under the group of isometries of $\Mm_k^d$. This group is the Gallilean group for the flat case, the orthogonal group $O(d+1)$ for the spherical case, and the Lorentz group $O(d,1)$ for the hyperbolic case. 

Christ \cite{Christ0} showed that $\RR_\nu$ is continuous from $L^p(\Mm_\nu^d)$ to $L^q(\MM_k(\Mm_\nu^d))$, where
\begin{equation}\label{eq:2a}
p \de \dfrac{d+1}{k+1}, \ \ \ q \de d+1.
\end{equation}
The associated sharp inequalities -- whose group of symmetries is the affine group $\Aff(\R^d)$ -- was studied in Christ \cite{Christ2}, Drouot \cite{Drouot1} and Flock \cite{Flock}. The resulting theorem is:

\begin{theorem}\label{thm:1}\cite{Christ2, Drouot1, Flock} For $p, q$ given by \eqref{eq:2a}, $0 \neq f \in L^p(\R^d)$,
\begin{equation}\label{eq:2c}
\dfrac{|\RR_0 f|_q}{|f|_p} \leq A_0 ,  \ \ A_0 = \dfrac{|\RR_0 h_0|_q}{|h_0|_p}, \ \ 
h_0(x) \de (1+|x|^2)^{-\frac{k+1}{2}}.
\end{equation}
Equality in \eqref{eq:2c} is realized if and only if
\begin{equation*}
f \in \{x \in \R^d \mapsto C h_0(Lx), \ L \in \Aff(\R^d), \ C \in \R\}.
\end{equation*}
\end{theorem}

Christ \cite{Christ2} first showed Theorem \ref{thm:1} when $k=d-1$. We derived later the sharp constant for any values of $k$ in \cite{Drouot1}. Flock \cite{Flock} completed the characterization by showing the uniqueness part of the statement.

Here we extend Theorem \ref{thm:1} to the spherical and hyperbolic $k$-plane transform. 

\begin{theorem}\label{thm:2} For $p,q$ given by \eqref{eq:2a}, $0 \neq f \in L^p(\Ss^d)$, 
\begin{equation*}
\dfrac{|\RR_+ f|_q}{|f|_p} \leq A_+ ,  \ \ A_+ = \dfrac{|\RR_+ \1_{\Ss^d}|_q}{|\1_{\Ss^d}|_p},
\end{equation*}
with equality realized if and only if
\begin{equation*}
f \in \left\{ \w = (\w',\w_{d+1}) \in \Ss^d \mapsto C \left( |\w_{d+1}|^2+|L\w'|^2 \right)^{-\frac{k+1}{2}}, \ L \in \Aff(\R^d), \ C \in \R \right\}.
\end{equation*}
\end{theorem}

\begin{theorem}\label{thm:3} For $p,q$ given by \eqref{eq:2a}, $0 \neq f \in L^p(\Hh^d)$,
\begin{equation}\label{eq:2b}
\dfrac{|\RR_- f|_q}{|f|_p} \leq A_- ,  \ \ A_- = \lim_{\lambda \rightarrow + \infty} \dfrac{|\RR_- h_\lambda|_q}{|h_\lambda|_p}, \ \ h_\lambda(\zeta) = \lambda^{d/p} \left(|\zeta_{d+1}|^2+\lambda^2 |\zeta'|^2\right)^{-\frac{k+1}{2}}.
\end{equation}
Equality in \eqref{eq:2b} is never realized. 
\end{theorem}

The proofs of Theorems \ref{thm:2} and \ref{thm:3} combine ideas of Drury \cite{Drury2} with the work on sharp forms of the $k$-plane transform on $\R^d$. Drury's work contains a very enlightening interpretation on the correspondance between the $k$-plane transform on the plane and on the sphere. His ideas can be transfered to the case of the hyperbolic plane. It is then at first a bit surprising that extremizers do not exist in this case, but this phenomena is standard while considering similar inequalities on $\Hh^d$ -- see Liu \cite{Li16} for work on the Sobolev inequality, and Banica--Duyckaerts \cite{Banica} for connections and blow-up results in the context of non-linear Schr\"odinger equations on $\Hh^d$.

Recently Chen \cite{Chen} showed the following multilinear inequality on $\Mm^d_\nu$: for any $\nu \in \{0,\pm\}$, $\tp \in (0,\infty)$, there exists $C > 0$ such that
\begin{equation}\label{eq:2g}
f_j \in L^{\tp}(\Mm_\nu^d) \ \ \Rightarrow \  \
\prod_{j=1}^{d+1} |f_j|_{L^{\tp}(\Mm_\nu^d)} \leq C \sup_{\zeta_j \in \Mm_\nu^d} \ \left(\prod_{j=1}^{d+1} f_j(\zeta_j)\right) \cdot |\det(\zeta_1, ...,  \zeta_{d+1})|^{\frac{d+1}{\tp}}.
\end{equation}
This inequality is closely related to \eqref{eq:2c} thanks to Drury's identity \eqref{eq:1u}.
Chen derived the best constant in \eqref{eq:2g} for $\R^d$ and $\Ss^d$, together with the value of some extremizers. The sharp form of \eqref{eq:2g} in the case of $\Hh^d$ remains untreated. In accordance with Theorem \ref{thm:3}, we conjecture that \eqref{eq:2g} on $\Hh^d$ admits no extremizers, and that a limit similar to \eqref{eq:2b} yields the best constant. The lacking ingredient is a full characterization of extremizing functions for Chen's inequality on $\R^d$, or at least a (rather weak) inverse result showing that if $f$ extremizes \eqref{eq:2g} on $\R^d$, then $f$ has full support. This remark will become clearer to the reader after looking at the proof of Theorem \ref{thm:3}. 

Every other inequality satisfied by $\RR_\nu$ on unweighted Lebesgue spaces can be obtained by interpolating the trivial $L^1 \rightarrow L^1$ bound with the $L^p \rightarrow L^q$ bound, where $p,q$ are given by \eqref{eq:2a}. Baernstein and Loss \cite{BaLo} made some conjectures about their sharp forms, proved in \cite{Christ2, Drouot1} in the endpoint cases. Proving them for all exponents would be a spectacular achievement -- similar questions remain open in the simpler case of the Hardy--Littlewood--Sobolev inequality. In the case of the Radon transform ($k=d-1$), Christ \cite{Christ0, Christ2} showed a stable form of \eqref{eq:2a}. This was extended to the $k$-plane transform in \cite{Drouot2}, with the extra assumption that the function $f$ is radial. Removing this restriction is an ambitious project, because of the large size of the group of symmetries $\Aff(\R^d)$. The first step should be an analysis of sets which have an overall large incidence with $k$-planes, in the spirit of Christ \cite{Christ3}.

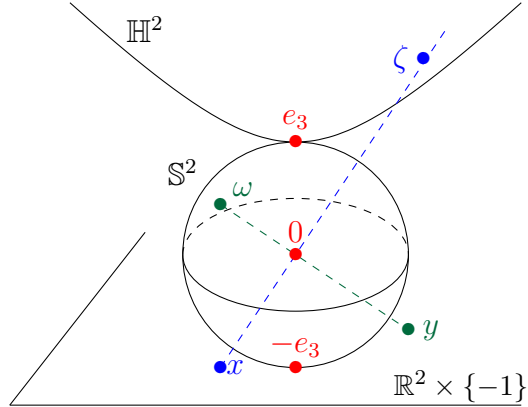
\begin{figure}
\centering
\begin{tikzpicture}
\draw (-2.8,-2) -- (4,-2);
    \draw (-2.8,-2) -- (-1,0.3);
\node[red] at (1,0.3) {$0$};
\node[red] at (1,-1.2) {$-e_3$};
    \draw (-.5,0) arc (180:360:1.5cm and 0.75cm);
    \draw[dashed] (-.5,0) arc (180:0:1.5cm and 0.75cm);
    \draw (1,0) circle (1.5cm);
\node at (-.5,1.1) {$\mathbb{S}^2$};
\node at (3.2,-1.75) {{\small $\R^2 \times \{-1\}$}};
\node[red] at (1,1.8) {$e_3$};
      \draw[domain=-2:4,smooth,variable=\t,black]  plot ({\t},{sqrt((\t-1)^2+1.5^2)});
\node at (-1,3) {$\mathbb{H}^2$};
\node[blue] at (.2,-1.5) {$x$};
\draw[blue, dashed] (0,-1.5) -- (3,3);
\draw[cadmiumgreen, dashed] (0,.67) -- (2.5,-1);
\node[cadmiumgreen] at (2.8,-1) {$y$};
\node[cadmiumgreen] at (.3,.87) {$\w$};
\node[blue] at (2.4,2.6) {$\zeta$};
\node[blue] at (2.7,2.6) {$\bullet$};
\node[blue] at (0,-1.5) {$\bullet$};
\node[cadmiumgreen] at (2.5,-1) {$\bullet$};
\node[cadmiumgreen] at (0,.67) {$\bullet$};
\node[red] at (1,1.5) {$\bullet$};
\node[red] at (1,-1.5) {$\bullet$};
\node[red] at (1,0) {$\bullet$ };
\end{tikzpicture}
\caption{ The spaces $\Ss^2, \ \Hh^2$ and $\R^2 \equiv \R^2 \times \{-1\}$ seen as subspaces of $\R^3$. In blue (resp. green) we plot $x = p_-(\zeta)$ for $\zeta \in \Hh^2$ (resp. $y = p_+(\w)$ for $\w \in \Ss^2$), which is the intersection of the line containing $\zeta$ (resp. $\w$) and $0$ with the plane $\R^2 \times \{-1\}$. Hence, $p_-$ maps great hyperbola in $\Hh^2$ to lines in $\R^2$, and $p_+$ maps great circles in $\Ss^2$ to lines in $\R^2$.}
\label{fig:31}
\end{figure}

Theorems \ref{thm:2} and \ref{thm:3} will follow from Proposition \ref{prop:1} below.  Define the projection $p_\pm :  \Mm_\pm^d  \rightarrow  \R^d$ by $p_\pm(\zeta',\zeta_{d+1}) = -\zeta'/\zeta_{d+1}$ -- a geometric interpretation is provided in Figure \ref{fig:31}. The pullback operator $p_\pm^*$ induces an operator $P_\pm$ given by
\begin{equation}\label{eq:1}
P_\pm f(\zeta) = p^*_\pm (\lr{x}^{k+1}_\pm f)(\zeta), \ \ \ \ f \in L^p(\R^d),
\end{equation}
where $\lr{x}_\pm = \sqrt{1\pm|x|^2}$.

\begin{proposition}\label{prop:1} After possibly multiplying the measures on $\Mm_\pm^d$ by a constant,
 \begin{itemize}
\item[$(i)$] For any $f \in L^p(\R^d)$,
\begin{equation}\label{eq:1c}
|P_+ f|_p = |f|_p,  \ \ |\RR_+ P_+ f|_q = |\RR_0 f|_q.
\end{equation}
\item[$(ii)$] For any $f \in L^p(\R^d)$ with support contained in the unit ball $\Bb^d \subset \R^d$,
\begin{equation}\label{eq:1e}
|P_- f|_p = |f|_p,  \ \ |\RR_- P_- f|_q = |\RR_0 f|_q.
\end{equation}
\end{itemize} 
\end{proposition}

The paper is organized as follows. Proposition  \ref{prop:1} shows that $P_\pm$ intertwines the operators $\RR_\pm$ with $\RR_0$. Theorem \ref{thm:2} is then a direct consequence of Proposition \ref{prop:1}, and of the bijectivity of $P_+$; we omit the proof. In contrast with the spherical case, $P_-$ does not act on the full $L^p(\R^d)$. This will generate the non-existence result of Theorem \ref{thm:3}, shown in \S\ref{sec:2}, where we assume that Proposition \ref{prop:1} holds. The proof of \eqref{eq:1e} in Proposition \ref{prop:1} is postponed to \S\ref{sec:3}. It relies on group-theoretic generators arguments and explicit computations. The proof of \eqref{eq:1c} is very similar and is omitted.

\noindent \textbf{Aknowledgement.} We thank M. Christ for bringing Drury's paper \cite{Drury2} to our attention. This research was supported by the NSF grant DMS-1500852 and the Fondation CFM pour la recherche.

\section{Proof of Theorem \ref{thm:3}.}\label{sec:2}

We prove Theorem \ref{thm:3} assuming that Proposition \ref{prop:1} holds. In this section $\MM^d_- = \Hh^d$.

 Let $P_-$ be the operator given by \eqref{eq:1}. Without loss of generalites, we can multiply the measure on $\Hh^d$ by a multiplicative constant so that $P_- : L^p(\R^d) \rightarrow L^p(\Hh^d)$ satisfies \eqref{eq:1e}. Since $p_-$ maps bijectively $\Bb^d$ to $\Hh^d$, the restriction of $P_-$ to functions with support in the unit ball is a bijection. Hence, the norm of $\RR_-$ must be at most equal to the norm of $\RR_0$. For $\lambda >0$, let $h_\lambda$ be the extremizer of \eqref{eq:2c} given by
\begin{equation*}
h_\lambda(x) = \lambda^{d/p} h_0(\lambda x), \ \ h_0(x) = \lr{x}_+^{-k-1}.
\end{equation*}
It is clear that $|h_\lambda|_p$ does not depend on $\lambda$ and that
\begin{equation*}
\lim_{\lambda \rightarrow +\infty} |h_\lambda|_{L^p(\R^d \mi \Bb^d)} = 0.
\end{equation*}
Therefore, if $\1_{\Bb^d}$ denotes the characteristic function of the unit ball, then
\begin{equation*}
A_- \geq \lim_{\lambda \rightarrow + \infty} \dfrac{|\RR_- P_- \1_{\Bb^d} h_\lambda|_q}{|P_- \1_{\Bb^d} h_\lambda|_p} = \lim_{\lambda \rightarrow + \infty} \dfrac{|\RR_0 \1_{\Bb^d} h_\lambda|_q}{|\1_{\Bb^d} h_\lambda|_p} \geq \lim_{\lambda \rightarrow + \infty} \dfrac{|\RR_0 h_\lambda|_q - A_0 |h_\lambda|_{L^p(\R^d \mi \Bb^d)}}{|h_\lambda|_p + |h_\lambda|_{L^p(\R^d \mi \Bb^d)}} = A_0.
\end{equation*}
This shows $A_- = A_0$ and \eqref{eq:2b}.

Assume now that an extremizer $u$ exists for \eqref{eq:2b}. Let $f \in L^p(\R^d)$ with support in the unit ball be such that $P_-f=u$. In particular, $|P_-f|_p = |u|_p$ and $A = |\RR_- u|_q = |\RR_- P_- f|_q = |\RR_0f|_q$: $f$ is an extremizer for $\RR_0$, with support in the unit ball. Christ--Flock's uniqueness result recalled in Theorem \ref{thm:1} shows that $f$ cannot be compactly supported, which is a contradiction.

\section{Proof of Proposition \ref{prop:1}.}\label{sec:3}

As said earlier, we prove only \eqref{eq:1e} in Proposition \ref{prop:1} -- which concerns the hyperbolic case $\nu = -$. Hence, in this section, $\Mm_\nu^d = \Hh^d$ and we drop all the subscripts ``minus" to make the notations simpler. 

\subsection{Preliminaries.}\label{sub:3.1}
We first observe that $p$ induces a bijection $\Hh^d \rightarrow \Bb^d$ and we define $p^{-1} : \Bb^d \rightarrow \Hh^d$ its inverse. For $\pi$ a $k$-plane in $\R^d$, $d(\pi)$ denotes the geodesic distance between $\pi$ and the origin of $\R^d$. If $d(\pi) < 1$, the set $p^{-1}(\pi)$ is a $k$-plane in $\Hh^d$, as the intersection of $\Hh^d$ with a $k+1$-plane in $\R^{d+1}$ -- see Figure \ref{fig:31}. For $f \in C^\infty(\R^d)$ with support in the unit ball, \cite[Theorem $3.1$]{BerTarKur} states that
\begin{equation} \label{eq:2d}
(\RRR P f) (p^{-1}(\pi)) = \lr{\pi} \RRR_0 f (\pi), \ \ \lr{\pi} \de \sqrt{1-d(\pi)^2}.
\end{equation}

If $O$ is an element of the Lorentz group $O(d,1)$, the composition operator with $O$ commutes with $\RR$: for $\xi \in \MM_k(\Hh^d)$, $u \in C^\infty_0(\Hh^d)$,
\begin{equation} \label{eq:2e}
(\RR u)(O \xi) = \RR (u \circ O) (\xi).
\end{equation}
The map $p$ conjugates the group $O(d,1)$ (acting naturally on elements of $\R^d$) to the group $\GG \de p O(d,1) p^{-1}$ (which acts naturally on elements of $\Bb^d$). In particular if $\Phi = p O p^{-1} \in \GG$, equations \eqref{eq:2d} and \eqref{eq:2e} imply that
\begin{equation}\label{eq:1y}
\lr{\Phi(\pi)} \RR_0 f (\Phi(\pi)) =\RR P f (O p^{-1}(\pi)) 
   = (\RR P (f \circ O)) (p^{-1}(\pi)) = \lr{\pi} \RR_0 \SSS f(\pi).
\end{equation}
In the above, $\SSS \de P O^* P^{-1}$ where $O^*$ denotes the pull-back operator with $O$: $\SSS$ is the $L^p$-isometry induced by the pull-back operator with $p O p^{-1}$.

Let $d\mu_{\MM_k(\Hh^d)}$ be a non-zero measure on $\MM_k(\Hh^d)$ that is invariant under $O(d,1)$. If $f : \R^d \rightarrow \R$ is a smooth function with support in the unit ball, the substitution $\xi = p^{-1}(\pi)$ shows that
\begin{equation*}
|\RR P f|_q^q = \int_{\MM_k(\Hh^d)} |\RR P f(\xi)|^q d\mu_{\MM_k(\Hh^d)}(\xi) = \int_{\MM_k(\R^d)} |\RR P f(p^{-1}(\pi))|^q d\mu_{\MM_k(\Hh^d)}(p^{-1}(\pi)). 
\end{equation*}
Since $f$ has support in $\Bb^d$, $\RR_0f(\pi) = 0$ when $d(\pi) \geq 1$. This together with \eqref{eq:2e} imply
\begin{equation}\label{eq:1o}
|\RR P f|_q^q = \int_{\MM_k(\R^d)} \lr{\pi}^q |\RR_0 f(\pi)|^q d\mu_{\MM_k(\Hh^d)}(p^{-1}(\pi)).
\end{equation}
In order to prove Proposition \ref{prop:1}, it suffices to show that the right hand side of \eqref{eq:1o} equals $|\RR_0 f|_q^q$, modulo a  multiplicative constant independent of $f$. From the point of view of measure theory, it is enough to show that $(p^{-1})^*d\mu_{\MM_k(\Hh^d)} = \lr{\pi}^{-q} d\mu_{\MM_k(\R^d)}$, where $(p^{-1})^*$ is the pull-back by $p^{-1}$. Since $d\mu_{\MM_k(\Hh^d)}$ is the unique measure on $\MM_k(\Hh^d)$ that is invariant under $O(d,1)$, $(p^{-1})^* d\mu_{\MM_k(\Hh^d)}$ is the unique measure on $\EE \de \{ \pi \in \MM_k(\R^d), d(\pi) \leq 1 \}$ that is invariant under $\GG$ -- modulo multiplicative constants. Therefore, to prove Proposition \ref{prop:1}, it suffices to show that $\lr{\pi}^{-q} d\mu_{\MM_k(\R^d)}$ is invariant under $\GG$; or equivalently, that for every $G \in C^\infty_0( \MM_k(\R^d))$ with $\supp(G) \subset \EE$, for every $\Phi \in \GG$,
\begin{equations}\label{eq:1f}
\int_{\MM_k(\R^d)} G(\pi) \lr{\pi}^{-q} d\mu_{\MM_k(\R^d)}(\pi) = \int_{\MM_k(\R^d)} G(\Phi(\pi)) \lr{\pi}^{-q} d\mu_{\MM_k(\R^d)}(\pi).
\end{equations}

We define the group $\OO(d)$ of $(d+1) \times (d+1)$ matrices of the form
\begin{equation}\label{eq:2}
\matrice{\Omega  & 0 \\ 0 & 1}, \ \ \Omega \in O(d),
\end{equation}
and the group $\BB$ of $(d+1) \times (d+1)$ matrices of the form
\begin{equation}\label{eq:3}
\matrice{\Id_{d-1} & 0 & 0 \\ 0 & a & -b \\ 0 & -b & a}, \  a^2-b^2=1,
\end{equation}
sometimes called boosts. Both $\OO(d)$ and $\BB$ are subgroups of $O(d,1)$ and we prove in the appendix that they generate $O(d,1)$ -- a result that must be known to specialists. Hence, it suffices to prove that \eqref{eq:1f} holds when $\Phi \in p \OO(d) p^{-1}$ and when $\Phi \in p \BB p^{-1}$. Since $p \OO(d) p^{-1} = O(d)$, \eqref{eq:1f} immediately holds for $\Phi \in p \OO(d) p^{-1}$. It remains to check \eqref{eq:1f} for transformations $\Phi \in p \BB p^{-1}$. We will verify it by performing explicit computations relying on Drury's identity \cite[Lemma 1]{Drury1}. This formula gives a fairly explicit description of $d\mu_{\MM_k(\R^d)}$: after possibly multiplying $d\mu_{\MM_k(\R^d)}$ by a multiplcative constant, for every $f \in C_0^\infty(\R^d)$, $F \in C_0^\infty(\MM_k(\R^d))$,
\begin{equations}\label{eq:1u} \\
\int_{\MM_k(\R^d)} |\RR_0 f(\pi)|^{k+1} F(\pi)  d\mu_{\MM_k(\R^d)}(\pi)
= \int_{\left(\R^d\right)^k} \dfrac{f(x_0) ... f(x_k)  F(\pi)}{\Det(x_0, ..., x_k)^{d-k}}dx_0 ... dx_k.
\end{equations}
Here, $\Det(x_0, ..., x_k)$ denotes the $k$-volume of the $k$-simplex with vertices $x_0, ..., x_k$ in $\R^d$, and $\pi$ is the $k$-plane containing $x_0, ..., x_k$ (this definition makes sense for almost every $x_0, ..., x_k$).

\subsection{Jacobian computations}

We fix here $\Phi \in p \BB p^{-1}$, and we denote by $J \Phi(x)$ the Jacobian of $\Phi$ at $x \in \R^d$. We show the following lemma:

\begin{lem}\label{lem:1}\begin{enumerate}
\item[$(i)$] There exist $a,b \in \R$ with $a^2 - b^2=1$ such that
\begin{equation}\label{eq:1r}
\Phi(x) = \left( \dfrac{ x'}{ bx_d+a}, \dfrac{ ax_d+b}{bx_d+a} \right), \ \ \ \ 
J\Phi(x) = \dfrac{1}{|bx_d+a|^{d+1}}.
\end{equation}
\item[$(ii)$] If $\SSS$ is the operator of \eqref{eq:1y}, then
\begin{equation}\label{eq:1s}
\SSS f(x) = \dfrac{1}{|bx_d+a|^{k+1}} f \circ \Phi(x).
\end{equation} 
\item[$(iii)$] If $x_0, ..., x_k \in \R^d$ are contained in a unique $k$-plane $\pi$,
\begin{equation}\label{eq:1g}
\dfrac{\lr{\pi}}{\lr{\Phi(\pi)}} = \dfrac{|bx_{0d}+a| \cdot ... \cdot |bx_{kd}+a| \cdot \Det(\Phi(x_0), ..., \Phi(x_k))}{\Det(x_0, ..., x_k)}.
\end{equation}
\end{enumerate}
\end{lem}

\begin{proof} For $(i)$: we have
\begin{equation*}
\Phi(x) = p \matrice{ \Id_{d-1} & 0 & 0 \\ 0 & a & -b \\ 0 &  -b & a } \dfrac{(-x,1)}{\lr{x}} = p \dfrac{(-x',-ax_d-b,bx_d+a)}{\lr{x}} = \left( \dfrac{ x'}{ bx_d+a}, \dfrac{ ax_d+b}{bx_d+a} \right).
\end{equation*}
Therefore
\begin{equation*}
J\Phi(x) =
\left| \begin{matrix}
\dfrac{1}{bx_d+a}\Id_{d-1} & 0 \\ \star & \dfrac{d}{dx_d}  \left(\dfrac{ax_d+b}{bx_d+a}\right)
\end{matrix}  \right| = \left|\dfrac{1}{(bx_d+a)^{d-1}} \dfrac{a^2-b^2}{(bx_d+a)^2}\right| = \dfrac{1}{|bx_d+a|^{d+1}}.
\end{equation*}

The statement $(ii)$ is a simple consequence of the fact that $\SSS$ is the $L^p$-isometry induced by the pull-back $\Phi^*$ and the form of $J\Phi(x)$ given by $(i)$.

The proof of $(iii)$ relies on ideas of \cite{Drouot1}; we believe that there is a more geometric proof. For $f \in C_0^\infty(\R^d)$, we define $\tRR_0 f : (\R^d)^{k+1} \rightarrow \R$ as
\begin{equation*}
\tRR_0 f(x_0, ..., x_k) = \int_{\lambda \in \R^k} f \left(x_0+\lambda_0(x_1-x_0) + ... + \lambda_k(x_k-x_0)\right) d\lambda.
\end{equation*} 
This function is related to $\RR_0 f$ by the formula 
\begin{equation*}
\Det(x_0, ..., x_k)\tRR_0 f(x_0, ..., x_k) = \RR_0 f (\pi),
\end{equation*}
see \cite[(3.5)]{Drouot1}. The identity \eqref{eq:1y} 
implies then
\begin{equation}\label{eq:1v}
\dfrac{\lr{\pi}}{\lr{\Phi(\pi)}} = \dfrac{\Det(\Phi(x_0), ..., \Phi(x_k))}{\Det(x_0, ..., x_k)} \cdot \dfrac{\tRR_0 f (\Phi(x_0), ..., \Phi(x_k))}{\tRR_0 \SSS f(x_0, ..., x_k)}.
\end{equation}
Therefore \eqref{eq:1g} is a consequence of a formula relating the quantities $\tRR_0 \SSS f(x_0, ..., x_k)$ and $\tRR_0 f (\Phi(x_0), ..., \Phi(x_k))$, that we shall prove below.

We start by writing $\delta x_k = x_k-x_0$ and $\az = x_0+\lambda_1 \delta x_1 + ... + \lambda_k \delta x_k$, so that
\begin{equation}\label{eq:1m}
\Phi(x_0+\lambda_1\delta x_1 + ... + \lambda_k \delta x_k) = \dfrac{\az'+(a\az_d+b)e_d }{b\az_d+a}.
\end{equation}
It implies that
\begin{equations}\label{eq:1x}
\tRR_0 \SSS f(x_0,...,x_k)  =  \int_{\R^k} \SSS f (x_0+\lambda_1\delta x_1 + ... + \lambda_k \delta x_k) d\lambda \\
  =  \int_{\R^k} f \left( \dfrac{\az'+(a\az_d+b)e_d }{b\az_d+a} \right) \dfrac{d\lambda}{|a\az_d+b|^{k+1}}.
\end{equations}
We make a first substitution $\lambda \mapsto \tlambda$ such that the argument $\frac{\az'+(a\az_d+b)e_d}{b\az_d+a}$ of $f$ in \eqref{eq:1x} is~linear~in~$\tlambda$:
\begin{equation*}
\tlambda_j = \dfrac{\tlambda_j}{b\az_d+a}- \dfrac{ax_{0d}+b}{bx_{0d}+a}, \ 1 \leq j \leq k-1, \ \ \  \lambda_k = \dfrac{a\az_d+b}{b\az_d+a}- \dfrac{ax_{0d}+b}{bx_{0d}+a}.
\end{equation*}
A computation using that $a^2-b^2=1$ shows that the argument $\frac{\az'+(a\az_d+b)e_d}{b\az_d+a}$ of $f$ in \eqref{eq:1x} depends indeed linearly on $\tlambda$:
\begin{equation}\label{eq:1a}
\dfrac{\az'+(a\az_d+b)e_d }{b\az_d+a} = \Phi(x_0) + \sum_{j=1}^{k-1} \tlambda_j \left( \delta x_j' + \dfrac{\delta x_{jd}}{\delta x_{kd}} \delta x_k' \right) + \tlambda_k \left( e_d + \dfrac{bx_{0d}+a}{\delta x_{kd}} \delta x_k' -b x'_0 \right).
\end{equation}

When $\lambda$ spans $\R^k$, $x_0+\lambda_1\delta x_1 + ... + \lambda_k \delta x_k$ spans $\pi$ and $\Phi(x_0+\lambda_1\delta x_1 + ... + \lambda_k \delta x_k)$ spans $\Phi(\pi)$. We deduce from \eqref{eq:1m} and \eqref{eq:1m} that $\frac{\az'+(a\az_d+b)e_d }{b\az_d+a}$ linearly spans $\Phi(\pi)$ when $\tlambda$ spans $\R^k$. It is then possible to make a linear substitution $\tlambda \mapsto \mu$ so that \eqref{eq:1a} becomes
\begin{equation}\label{eq:1n}
\dfrac{\az'+(a\az_d+b)e_d }{b\az_d+a} = \Phi(x_0) + \sum_{j=1}^k \mu_j(\Phi(x_j)-\Phi(x_0)).
\end{equation}
To this end, we define $A$ the $k \times k$ matrix such that $Af_i = e_i$, where $e_i$ is the canonical basis of $\R^k$ and 
\begin{equation*}
\begin{matrix}
f_i \de \dfrac{1}{bx_{id}+a} e_i + \left(\dfrac{ax_{id}+b}{bx_{id}+a} - \dfrac{ax_{0d}+b}{b x_{0d}+a} \right)e_k,\ \ \ \ 1 \leq i \leq k-1, \\
f_k \de \left(  \dfrac{ax_{kd}+b}{bx_{kd}+a} - \dfrac{ax_{0d}+b}{b x_{0d}+a} \right)e_k.
\end{matrix}
\end{equation*}
The substitution $\mu = A\tlambda$ and a computation yields \eqref{eq:1n}. The combination of \eqref{eq:1x} with \eqref{eq:1n} shows that
\begin{equation}\label{eq:1w}
\tRR_0 \SSS f(x_0,...,x_k) = \int_{\R^k} f\left(\Phi(x_0) + \sum_{j=1}^k \mu_j(\Phi(x_j)-\Phi(x_0))\right) 
\dfrac{|\det A^{-1}|}{\left| \det\partial \tlambda/ \partial \lambda \right|} \dfrac{d\mu}{|a\az_d+b|^{k+1}}.
\end{equation}

To conclude we need to compute the Jacobian determinants $\det A^{-1}$ and $\det\partial \tlambda/ \partial\lambda$.
Since $A^{-1}$ is upper triangular, its determinant is the product of its diagonal elements:
\begin{equation}\label{eq:10}
\det A^{-1} = \dfrac{x_{kd}-x_{0d}}{(bx_{0d}+a) \cdot ... \cdot (bx_{kd}+a)}.
\end{equation}
The same method as in \cite[Appendix]{Drouot1} leads to
\begin{equation}\label{eq:1p}
\det\left(\dd{\tlambda}{\lambda}\right) = \dfrac{x_{kd}-x_{0d}}{(b\az_d+a)^{k+1}}.
\end{equation}
We plug \eqref{eq:10} and \eqref{eq:1p} in \eqref{eq:1w} to obtain
\begin{equation*}
\tRR_0 \SSS f\left(x_0, ..., x_k\right) = \dfrac{\tRR_0 f \left(\Phi(x_0), ..., \Phi(x_k) \right)}{|bx_{0d}+a|\cdot...\cdot|bx_{kd}+a|}.
\end{equation*}
This relation is now plugged in \eqref{eq:1v}, showing \eqref{eq:1g}.\end{proof}

\subsection{Proof of Proposition \ref{prop:1}.}

\begin{proof}[Proof of Proposition \ref{prop:1}] Because of the discussion in \S \ref{sub:3.1}, we fix $G \in C^\infty_0(\MM_k(\R^d))$ with $\supp(G) \subset \EE$ and $\Phi \in p \BB p^{-1}$, and we show that \eqref{eq:1f} holds. Let $f$ be a smooth radial nonnegative function on $\R^d$, that is positive inside the (open) unit ball and vanishes outside the (open) unit ball. The function $\RR_0 f$ has support in $\EE$ and never vanishes on $\EE$. This allows to define $F \de |\RR_0 f|^{-k-1} G$. Below we will always denote $\pi$ the $k$-plane containing $x_0, ... x_k$. Applying successively \eqref{eq:1y}, Drury's identity \eqref{eq:1u}, and \eqref{eq:1s}, we get
\begin{equations*}
\int_{\MM_k(\R^d)} G(\Phi(\pi)) \lr{\pi}^{-q} d\mu_{\MM_k(\R^d)}(\pi) = \int_{\MM_k(\R^d)} |\RR_0 f(\Phi(\pi))|^{k+1} F(\Phi(\pi)) \lr{\pi}^{-q} d\mu_{\MM_k(\R^d)}(\pi) \\
= \int_{\MM_k(\R^d)} |\RR_0  \SSS f(\pi)|^{k+1} \dfrac{F(\Phi(\pi)) \lr{\pi}^{k-d}}{\lr{\Phi(\pi)}^{k+1}}  d\mu_{\MM_k(\R^d)}(\pi) \\
= \int_{\left(\R^d\right)^k} \SSS  f(x_0) ... \SSS f(x_k)  \dfrac{F(\Phi(\pi)) \lr{\pi}^{k-d}}{\lr{\Phi(\pi)}^{k+1}}  \Det(x_0, ..., x_k)^{k-d}dx_0 ... dx_k \\
= \int_{\left(\R^d\right)^k}  f(\Phi(x_0)) ... f(\Phi(x_k))  \dfrac{F(\Phi(\pi)) \lr{\pi}^{k-d}}{\lr{\Phi(\pi)}^{k+1}}  \cdot  \dfrac{\Det(x_0, ..., x_k)^{k-d} dx_0 ... dx_k}{\left( |bx_{0d}+a| \cdot ... \cdot |bx_{kd}+a|\right)^{k+1}}.
\end{equations*}
We now apply \eqref{eq:1g} and we substitute $y_i = \Phi(x_i)$ with corresponding Jacobian factor given by \eqref{eq:1r} to obtain:
\begin{equations*}
\int_{\left(\R^d\right)^k}  f(\Phi(x_0)) ... f(\Phi(x_k))  \dfrac{F(\Phi(\pi)) \lr{\pi}^{k-d}}{\lr{\Phi(\pi)}^{k+1}} \cdot \dfrac{\Det(x_0, ..., x_k)^{k-d} dx_0 ... dx_k}{\left( |bx_{0d}+a| \cdot ... \cdot |bx_{kd}+a|\right)^{k+1}} \\
= \int_{\left(\R^d\right)^k} f(\Phi(x_0)) ... f(\Phi(x_k)) \dfrac{F(\Phi(\pi))}{\lr{\Phi(\pi)}^q} \cdot \dfrac{\Det(\Phi(x_0), ..., \Phi(x_k))^{k-d} dx_0 ... dx_k}{\left( |bx_{0d}+a| \cdot ...\cdot  |bx_{kd}+a|\right)^q} \\
= \int_{\left(\R^d\right)^k}  f(y_0) ... f(y_k)  \dfrac{F(\pi)}{\lr{\pi}^{-q}}  \Det(y_0, ..., y_k)^{k-d} dy_0 ... dy_k.
\end{equations*}
Another use of Drury's identity \eqref{eq:1u} finally leads to
\begin{equations*}
\int_{\MM_k(\R^d)} G(\Phi(\pi)) \lr{\pi}^{-q} d\mu_{\MM_k(\R^d)}(\pi) = \int_{\MM_k(\R^d)} |\RR_0 f(\Phi(\pi))|^{k+1} F(\Phi(\pi)) \lr{\pi}^{-q} d\mu_{\MM_k(\R^d)}(\pi) \\ = \int_{\MM_k(\R^d)} |\RR_0 f(\pi)|^{k+1} F(\pi) \lr{\pi}^{-q} d\mu_{\MM_k{\R^d}}(\pi) = \int_{\MM_k(\R^d)} G(\pi) \lr{\pi}^{-q} d\mu_{\MM_k(\R^d)}(\pi).
\end{equations*}
This proves \eqref{eq:1f} and therefore Proposition \ref{prop:1}.
\end{proof}

\section{Appendix: a generator Lemma}

Recall that $\OO(d)$ is the subgroup of $O(d,1)$ given by matrices of the form \eqref{eq:2} and $\BB$ is the subgroup of $O(d,1)$ given by matrices of the form \eqref{eq:3}.

\begin{lem} The group $O(d,1)$ is generated by $-\Id$, $\OO(d)$ and $\BB$.
\end{lem}

\begin{proof} According to the Dieudonn\'e--Cartan theorem \cite[pp 10]{Ca}, every element of $O(d,1)$ is the product of (at most $d+1$) reflections; we recall that reflections in $O(d,1)$ are transformations of the form
\begin{equation*}
R_n : v \mapsto v - 2 \dfrac{\lr{v,n}_-}{\lr{n,n}_-^2}n, \ \ \  n \in \R^{d+1}, \  \ \ \lr{x,y}_- \de x_1y_1 +... +x_dy_d-x_{d+1}y_{d+1}.
\end{equation*}
Hence, it suffices to prove that each reflection $R_n$ can be written as the products of elements in $\OO(d)$ and $\BB$. We observe that if $L \in O(d,1)$ then
\begin{equation}\label{eq:1z}
L^{-1} R_n L v = v - 2 \dfrac{\lr{Lv,n}_-}{\lr{n,n}_-^2} L^{-1} n = v - 2 \dfrac{\lr{v,L^{-1} n}_-}{\lr{L^{-1}n,L^{-1}n}_-^2} L^{-1} n = R_{L^{-1}n} v.
\end{equation}
Let $\Omega \in \OO(d)$ such that $\Omega n \in \{0\}^{d-1} \times \R^2$, and $B \in \BB$ so that $B\Omega n = \lr{n,n}_-e_{d+1}$. The composition formula \eqref{eq:1z} shows that $R_n = B^{-1} \Omega^{-1} R_{e_{d+1}} \Omega B$. In addition, $R_{e_{d+1}}$ is the product of $-\Id$ with $-\Id_d \oplus \Id_1 \in \OO(d)$; this shows that reflections are generated by $-\Id$, $\BB$ and $\OO(d)$. This completes the proof of the lemma.
\end{proof}

\end{document}